\documentclass[reqno]{gokova}

\input{goksty.aaa}
\usepackage[dvips]{graphicx,color}
\usepackage{amssymb}
\usepackage{latexsym}
\usepackage{euscript,epsfig}
\usepackage[arrow,matrix,line]{xy}

\newcommand{\zz}{\mathbb{Z}}
\newcommand{\rot}{\mathbb{T}}

\title{Symplectic surgeries from singularities}
\author[SMITH, THOMAS]{Ivan Smith and Richard Thomas}

\thanks{I.S. is supported by an EC 
    Marie-Curie fellowship No. HPMF-CT--2000-01013.  R.T.
    is supported by a Royal Society
    University Research Fellowship.}

\begin{document}

\address{Ecole Polytechnique, Palaiseau, France}
\email{smith@math.polytechnique.fr}

\address{Imperial College, London, UK}
\email{rpwt@ic.ac.uk}

\volume{9}

\maketitle

\section{The local model}

Given an isolated analytic hypersurface singularity $0 \in X_0:=
\{f(z)=0\} \subset \C^{n+1}$ one can
form the \emph{smoothing} $X_t=\{f(z)=t\}$ and the \emph{resolution} $\widehat
X\to X_0$, obtained by (repeatedly) blowing up the origin. These two associated
spaces have the same link -- that is, the intersection $\SS^{2n+1} \cap
f^{-1}(0)$ -- and there is a smooth surgery which replaces the
smoothing by the resolution, or vice-versa. Thinking of the smoothing as
K{\"a}hler (and so symplectic) by restriction of the standard K{\"a}hler form
on $\C^{n+1}$, the smoothing also looks somewhat like a resolution by
symplectic parallel transport, described below. The result is that there is a
canonical map or `Lagrangian blow up' $X_0\leftarrow X_t$ whose `exceptional
locus' (the inverse image of the singular point $0\in X_0$) is a
\emph{Lagrangian} cycle,
in fact a collection of Lagrangian spheres \cite{Se3}. So the surgery
replaces configurations of Lagrangian spheres by complex (symplectic)
subvarieties. 
Symplectic parallel transport also shows that
the $X_t$s are all isomorphic \emph{as symplectic manifolds}, so denoting
any such symplectic manifold by $X$, we can denote this surgery
by the diagram (motivated by smoothings and resolutions in algebraic geometry)
\begin{equation} \label{ivanisagirl}
\xymatrix@R=.85em@C=.85em{
\widehat X \ar[d] \\
X_0 & X. \ar[l]
}\end{equation}
More general singularities, including complete intersections and some
non-isolated singularities, can be treated similarly; for simplicity
we will restrict 
attention to isolated hypersurface singularities.

We briefly remind the reader about symplectic parallel transport \cite{Se3}.
The total space $\mathcal X=\{X_t\}_{t\in\C}\stackrel{p}{\longrightarrow}\C$
is naturally a smooth symplectic submanifold of $\C^{n+2}$:
$\mathcal X=\{(z,t)\in\C^{n+1}\times\C\,:\,f(z)=t\}\subset\C^{n+2}$.
Away from the origin, there is a natural connection on this family of
$X_t$s, whose horizontal 
subspaces of the tangent bundle $T\mathcal X$ are the annihilators under
the symplectic form of $Tp = \ker(dp)$, the fibrewise tangent bundles
of the fibres $TX_t$. 
Parallel transport along this connection identifies fibres $X_t\cong
X_s$ (for $s,t \neq 0$),
and preserves the symplectic form on the fibres, so that the isomorphism
is of symplectic manifolds. Following the transport as $t\to0$ gives in
the limit the
map $X_t\to X_0$ above, collapsing the Lagrangian sphere vanishing cycles
of the singularity. (That these cycles are spheres follows from Morsifying
$f$ before totally smoothing it to $f-t$, and invoking the situation 
for ordinary double points as described in \cite{Se3}.)

The resolution is also K{\"a}hler (and so symplectic) as the blow up of $X_0$
sits naturally in $\C^{n+1}\times\P^n$. (\emph{Small} resolutions may
not be naturally
K{\"a}hler, however, or there may be choices involved, leading to obstructions
in patching these choices in a global situation; more of this later.) So
in a natural way the surgery (\ref{ivanisagirl}) is symplectic. This can
also be seen in a different way, as an instance of ``gluing along
convex boundaries''.  Removing a small tubular neighbourhood of 
the Lagrangian vanishing cycles, the boundary of the resulting manifold
is $\omega$-convex \cite{Et} and (by parallel transport) contactomorphic to the
$\omega$-convex link of the
exceptional divisor in the blow-up.  The Lagrangian and complex
fillings of this link can then be exchanged by a symplectic surgery
\cite{Et}.  In many cases of interest (for
instance singularities arising from weighted homogeneous polynomials)
the link is fibred by circles -- leaves of the
characteristic foliation -- on
which the symplectic form is degenerate. Quotienting out by these circles
gives a space with a symplectic form; this is the blow up of the
singularity.  (If it is singular we can blow up again.) The choice of the
size of the neighbourhood of the vanishing cycles goes over to the
size of the symplectic form
on the $\P^n$ factor of $\C^{n+1}\times\P^n\supset\widehat X$.

This paper describes work still in progress;
we discuss several surgeries which fit into the above framework, but
the examples to which they give rise need further study.  
For now, we will motivate the idea that they should be useful
in addressing various existence questions for symplectic
structures. After clarifying the global symplectic geometry of the
surgeries below, the subsequent sections deal with surgeries
relevant to symplectic manifolds with $c_1>0$ (Fano), $c_1=0$ (Calabi-Yau)
and $c_1<0$ (general type) respectively.
\bigskip 

\section{The global model}

Suppose we have a symplectic manifold $X$ containing a
configuration $\mathcal{C}$ of Lagrangians,
with a neighbourhood $U(\mathcal{C})$ isomorphic to the neighbourhood
of the vanishing cycles
of a K{\"a}hler manifold which has a K{\"a}hler degeneration collapsing
$\mathcal{C}$ to an isolated hypersurface singularity. (We discuss
below conditions 
which can ensure this; see Proposition \ref{maslov} for instance.) More
generally, suppose $X$ contains finitely many disjoint
such configurations, all
locally the vanishing cycles of (possibly different) K{\"a}hler degenerations.

Collapsing
these Lagrangians gives a singular space with a neighbourhood of each
singular point $p$ isomorphic to a neighbourhood $V(p)\subset\C^{n+1}$
of the above singularity.  Moreover
this isomorphism takes the symplectic form $\omega_X$ on $X$ to the restriction
of $\omega_{\C^{n+1}}$ to $V(p)$.
Pulling back $\omega_X$ to the blow-up $\widehat X$ of $X$ at the
singular points gives a 2-form $\omega$ degenerate
only along the exceptional locus $E$.

Assume for simplicity this blow-up $\widehat X$ is smooth; if not we can
iterate the following process. Let $\sigma$ be a closed 2-form on $\widehat X$
Poincar{\'e} dual to $-[E]$ supported on the neighbourhood of $E$ that is
the pull back of the K{\"a}hler neighbourhood of the singularity above. On
this neighbourhood $\sigma$ is cohomologous
to the K{\"a}hler form used in the local model above (restricted
from $\C^{n+1}\times\P^n$), and so can be taken to be equal to this
K{\"a}hler form in a (smaller) neighbourhood $U$ of $E$. Then we claim that
for $0<\varepsilon\ll1,\ \omega+\varepsilon\sigma$ is symplectic
globally on $\widehat X$.

Since nondegeneracy is an open condition, this is clear on $\widehat
X\backslash U$ for small enough $\varepsilon$, and on $E$ itself. On
$U\backslash E$ it follows from the fact that both $\sigma$ and
$\omega$ are compatible with the local complex structure inherited
from the local K{\"a}hler model; therefore any convex linear combination
of the two forms is also symplectic by an observation of Gromov.
 
The size of $\varepsilon$ is determined by the areas of curves inside the
exceptional divisor $E$, which in turn is related to the volume of the
neighbourhood $U(\mathcal{C})$.  This phenomenon is well-known from
symplectic blowing up at smooth points \cite{McD-S},
cf. (\ref{floppyears}) below. \bigskip 

\section{The ordinary double point}

Let $(X,\omega)$ be a symplectic manifold and $L\subset X$ a
Lagrangian sphere.  According to a theorem of Weinstein, a
neighbourhood of $L$ in $X$ is symplectomorphic to a neighbourhood of
the zero-section $L_0$ in the cotangent bundle $T^*L$, equipped with its
canonical symplectic structure.  Slightly less well known is the
existence of a symplectomorphism
\begin{equation} \label{cotangentisnode}
\left( \left\{ \sum_{i=1}^{n+1} z_i^2 = 0 \right\} \backslash \{0\},
\frac{i}{2} \sum dz_j
\wedge d\overline{z}_j \right) \ \cong \  (T^*\SS^n \backslash
\SS^n_0, dp\wedge dq) 
\end{equation}
\noindent although this can be given in a straightforward manner in
co-ordinates, for instance by taking $(z_j=a_j+ib_j)_j \mapsto
(a_j/|a|, -|a|b_j)_j$.  Here we have used the round metric to identify
$T^*\SS^n$ and $T\SS^n$.  The same map defines a global isomorphism from
$T^*\SS^n$ to $\{\sum z_j^2 = t\}$ when $t$ is real and positive,
explicitly exhibiting the
cotangent bundle as a smoothing of the singularity (and there is a similar
map for all $t\in\C^*$).  The space $W$ on the left hand side of 
(\ref{cotangentisnode}) is a punctured neighbourhood of the
$n$-fold \emph{ordinary double point} or \emph{node}.  This admits a
holomorphic desingularisation by blowing up the origin, which gives an
exceptional divisor a complex quadric $Q_{n-1} \subset \P^n$.   The total
space of the normal bundle $\mathcal{L}$ to this quadric in the
resolution admits two natural maps: the projection to $Q_{n-1}$, and a map
to $\C^{n+1}$ (whose image lies inside the singular space $W$), which
we label
$$Q_{n-1} \ \stackrel{p}{\longleftarrow} \ \mathcal{L} \
\stackrel{\pi}{\longrightarrow} \ \C^{n+1}.$$
We can define a model family of symplectic forms on the neighbourhood 
of the exceptional divisor in the resolution by setting
$$
\rho_{\lambda} = \pi^* \omega_{\C^{n+1}} + \lambda^2 p^*
\omega_{Q_{n-1}}^{\ }
$$
where the form on $Q_{n-1}$ is the restriction of the Fubini-Study
form\footnote{Our 
  convention throughout the paper is that the (Fubini-Study) symplectic form on
  $\P^n$ is given by quotienting the Hopf circles in the unit sphere
  in $\C^{n+1}$ and descending the standard form $\sum dx_j \wedge
  dy_j$, and this gives a line $\P^1 \subset \P^n$ area $\pi$.}
to $Q_{n-1}\subset\P^n$.  The form
$\rho_{\lambda}$ gives the generators of $H_2(Q_{n-1})$ equal size
$\pi \lambda^2$ (the generator is unique unless $n=3$).  Moreover let
$B(\delta)$ denote the ball of radius $\delta$ in $\C^{n+1}$.

\begin{lem} \label{floppyears}
There is a symplectomorphism between $[\pi^{-1} B(\delta)\backslash Q_{n-1},
\rho_{\lambda}]$ and the ``shell'' $[(B(\sqrt{\lambda^2+\delta^2})\backslash
B(\lambda))\cap W, \omega_{\C^{n+1}})]$.
\end{lem}

\begin{proof}
There is a diagram $\P^n \leftarrow \tilde{\mathcal{L}} \rightarrow \C^{n+1}$
arising from blowing up the origin in $\C^{n+1}$, with $\tilde{\mathcal{L}}$
the total space of the $\mathcal{O}(-1)$ line bundle over $\P^n$.
Denote by $\tilde{\rho}_{\lambda}$ the form
$\tilde{\pi}^* \omega_{\C^{n+1}} + \lambda^2 \tilde{p}^*
\omega_{\P^n}$, in an obvious notation.  Then according to 
(\cite{McD-S}, Lemma 7.11 or Lemma 6.40 in the 1st edition) there is a
symplectomorphism $z\mapsto F\circ\tilde{\pi}(z)$ between
$(\tilde{\pi}^{-1} B(\delta) \backslash \P^n,
\tilde{\rho}_{\lambda})$ and $(B(\sqrt{\lambda^2+ \delta^2})\backslash
B(\lambda), \omega_{\C^{n+1}})$, where $F$ is the radial map $z \mapsto
(\sqrt{|z|^2+\lambda^2}/|z|) z$.  The map $F$ preserves the quadric
$W$, so there is an induced map 
$F\circ\pi$ between the spaces given in the Lemma.  Since the
symplectic structures on these are induced from the ambient $\P^n
\times \C^{n+1}$ by restriction, this is a symplectomorphism.
\end{proof}

\noindent This formalises the way in which we can perform the obvious
surgery, by removing $B(\sqrt{\lambda^2+\delta^2})$ and gluing back
$\pi^{-1}B(\delta)$, in a
manner compatible with symplectic forms.  It is worth emphasising that
although the blow-up of a singular projective variety will always
remain projective (hence K{\"a}hler), in general there may be more
symplectic degenerations of a projective variety than exist
holomorphically; we can perform the symplectic surgery above starting
with any Lagrangian sphere, and not just a vanishing cycle for a
complex degeneration.  Of course, the difference between these two
classes (if any) is largely mysterious.  \medskip
   
\noindent The formalism above motivates the question of
determining the \emph{maximum possible} value of $\lambda$
that one can take in (\ref{floppyears}).  This is analogous to ``symplectic
packing'' questions, only looking not for symplectically standard balls
but for symplectically standard (Lagrangian) disc bundles over Lagrangian
submanifolds.  One reason to focus on this variant of a packing number
is given by the following (itself a variant of ideas of symplectic
inflation).  Fix once and for all a model of
$T^*\SS^n$ as $\{(u,v) \in \R^{n+1}\times \R^{n+1} \, | \, |u|=1,
\, \langle u,v\rangle =0 \}$, and regard all Lagrangian spheres as
parametrised by the zero-section in this model equipped with the
standard form $du \wedge dv$.  Let us say a symplectic manifold is a
\emph{symplectic Fano} 
if $[\omega] = c_1(TX,\omega) \in H^2(X,\zz)$.  These are essentially
the well-known
``monotone'' symplectic manifolds, although to fix the constants in the
next Lemma it is important that
$[\omega]$ and $c_1$ co-incide and are not just positively proportional.

\begin{lem} \label{perhapsfano}
Let $(X,\omega)$ be a symplectic Fano 6-manifold and $L\subset X$ a Lagrangian
sphere.  If $\{v \in T^*\SS^3 \ | \ |v| \leq \mu \}$ 
symplectically embeds inside $X$ for any $\mu > \frac{1}{2\pi}$ then
the manifold $Y$ obtained by surgery along $L$ 
admits symplectic structures in the cohomology class $c_1(Y)$.
\end{lem}

\begin{proof}
The first Chern class of $X$ is zero on the Lagrangian $L$, so lifts to
a class $c_1(X)\in H^2(X,L)\cong H^2(Y,E)\to H^2(Y)$, with
$E\cong\P^1\times \P^1$ the exceptional divisor. By the adjunction formula
its image satisfies $c_1(Y) = c_1(X) - [E]$.
On the other hand, the cohomology class of the form
$\rho_{\lambda}$ given above is given by $[\omega_X] - \pi \lambda^2
E$.  (We again identify $[\omega_X]$ with an element of $H^2(Y)$, and
the factor of $\pi$ enters because a line in projective space,
equipped with the Fubini-Study form as per our conventions, has area $\pi$.)
We would therefore
like to take $\lambda = 1/\sqrt{\pi}$ in performing the surgery. On
the other hand, using the
local model (\ref{cotangentisnode}), one can check that the ball
$(B(\sqrt{1/\pi + \delta^2}), \omega_{\C^3})$ symplectically embeds
inside $\{v \in T^*\SS^3 \ | \ |v| \leq 1/2\pi + \delta^2/2 \}$.   This is just
because $|z|^2 \leq R \Rightarrow |v| = |\Re(z)||\Im(z)| \leq R/2$.
The result follows.
\end{proof}

\noindent  In four dimensions, it
is known that every symplectic Fano is in fact a del Pezzo surface
\cite{Liu}, and so K{\"a}hler. (The above surgery would not be relevant in
4 dimensions, as it produces a symplectic $-2$-curve on which $c_1$ is
zero, not positive.) In higher dimensions, there is no analogous result, nor
is any counterexample known.  Lemma (\ref{perhapsfano}) provides a symplectic
surgery that preserves the class of symplectic Fanos.  It can be applied
in two directions. On the one hand, there is a classification of Fano
3-folds \cite{MM} and one can look for symplectic non-K{\"a}hler Fanos;
on the other, restrictions on symplectic Fanos will translate
into packing-type bounds for neighbourhoods of Lagrangian spheres.
For instance, there is a Lagrangian sphere
inside $(\P^2
\times \P^1, \frac{1}{\pi}(-3\omega_{FS} \oplus 2\omega_{FS}))$ given as
follows.  Embed a ball $B(\sqrt{2/\pi}) \subset \P^2$ with the
standard Euclidean symplectic form on the left and $3/\pi$ times the
Fubini-Study form on the right.  The boundary $\SS^3$ of this ball
maps into $\P^2 \times \P^1$ via the graph of the Hopf map, and is
Lagrangian with respect to the chosen form (which induces the usual
orientation on each factor and is normalised so
that $c_1 = [\omega]$). Alternatively, we can remove the minus sign, yielding a
symplectic form on $\P^2 \times \P^1$ deformation equivalent to the
usual K{\"a}hler form, if we compose 
the Hopf map with the antipodal map of $\P^1$ in the definition of
the Lagrangian sphere.

\begin{ques}
Is there a symplectic embedding of the $\mu$-disc bundle of
$T^*\SS^3$ into $(\P^2
\times \P^1, \frac{1}{\pi}(-3\omega_{FS} \oplus 2\omega_{FS}))$ for
any $\mu > 1/2\pi$?
\end{ques}

\noindent If the answer to this question is yes, then there are
symplectic Fano manifolds which are not K{\"a}hler.  For according to
the previous Lemma, the
transition of $\P^2 \times \P^1$ in the Lagrangian would be a
symplectic Fano.  This manifold would have $b_2=4$, with two new $H_2$
classes coming from the rulings of the exceptional divisor (and
classes in $H_4$ coming from the divisor, and the lift of a 4-chain
bounded by the Lagrangian).  However, the almost complex structure underlying
this symplectic structure has $c_1^3 = 52$, which is not realised by
\emph{any} Fano 3-fold with $b_2=4$, according to the classification lists of
\cite{MM}\footnote{Recently, Mori has announced a gap in \cite{MM}; there is 
an additional Fano 3-fold with $b_2=4$, but with $c_1^3 = 26$, given by the
blow-up of $\P^1\times\P^1\times\P^1$ along a $(1,1,3)$ curve.}.

On the other hand, if the question has a negative answer,
this gives a new kind of packing
obstruction: for the  volume of the $(1/2\pi)$-disc bundle is
strictly less than the volume of $\P^2 \times \P^1$ with the Fano
symplectic form.  \medskip

\noindent Incidentally, the result of the surgery above again contains
a symplectic two-sphere with trivial normal bundle (coming from the
diagonal curve inside $\P^1\times\P^1$ viewed inside the total space
of $\mathcal{O}(-1,-1)$).  However, a few moment's reflection
with Gromov's non-squeezing theorem shows that the symplectically
trivial neighbourhood of this curve is not large enough to
(necessarily) contain a Lagrangian sphere by the prescription above,
so we cannot expect to iterate the surgery symplectically. \bigskip

\section{Small resolutions}

The ordinary double point has special features in six real dimensions,
where the exceptional divisor $Q_2$ has two dimensional second
homology.  In fact, $Q_2 \cong \P^1 \times \P^1$ and there are
\emph{small resolutions} of the 3-fold node in which either of the two
rulings of $Q_2$ are contracted (i.e. replace the singular point by
just a rational curve).  The two possible small resolutions differ by
a flop, as described extensively in \cite{STY}.  Explicitly, writing
the node as $\{xy=zw\} \subset \C^4$, the small resolutions are given
by taking the graphs of two distinct maps to $\P^1$, namely $(x/z =
w/y)$ and $(x/w = z/y)$.  

Now in contrast to
blow-ups, small resolutions are not operations within the projective
or K{\"a}hler category, but they have another special feature: they
preserve the first Chern class of the manifold. This
gives a route to searching not for exotic Fano manifolds, as above,
but exotic Calabi-Yaus, at least under appropriate conditions for the
surgery to exist symplectically at all.  These were investigated in
\cite{STY}, where the surgeries given by replacing a Lagrangian
3-sphere by a symplectic two-sphere were called ``conifold transitions''.

\begin{thm} \label{dickisgay} Let $X$ be a symplectic
  six-manifold containing disjoint Lagrangian 
  spheres $L_1, \ldots, L_n$ for which $\sum \lambda_i [L_i] = 0 \in
  H_3(X,\zz)$ (with all $\lambda_i \neq 0)$.  Then at least one of the
  conifold transitions of $X$ in the
  $\{L_i\}$ admits a symplectic structure.
\end{thm}

\begin{proof}[Sketch]
We regard the existence of the homology relation as given by a
four-chain with boundary the union of the
Lagrangian spheres.  The boundary is collapsed in the nodal space,
giving a closed four-cycle which lifts to any small resolution.
Flopping gives a cycle $D$ hitting each resolving $\P^1$ positively
(with intersection number $|\lambda_i|$).
This cycle now plays the role taken by the exceptional divisors of
blow-ups in the argument given in ``The Global Model''; it is at least
cohomologically positive on the exceptional curves. Some local analysis
then shows that we can choose a two-form $\sigma$ Poincar{\'e} dual
to $D$ such that the form $\omega_X + \varepsilon \sigma$ is globally
symplectic for all $0<\varepsilon\ll1$: see \cite{STY} for details.
\end{proof}

\noindent As explained in \cite{STY}, this theorem is related via
  mirror symmetry to an operation on complex 3-folds studied by
  Clemens, Friedman and Tian -- this operation being the obvious
  inverse process, in which rational curves are contracted to give
  ordinary double points, and one looks for sufficient conditions to
  have a complex smoothing of the resulting nodal variety.  It is not
  hard to find examples in which the theorem can be applied; for
  instance, the Lagrangian 3-sphere exhibited as the graph of a Hopf
  map, embedded in
  $X \times \P^1$ for any symplectic four-manifold $X\supset \SS^3$ and an
  appropriate product symplectic form $\omega_X \oplus
  (-\omega_{FS})$, certainly satisfies $[L]=0 \in H_3$.  In general,
  the resulting manifolds can be K{\"a}hler; for instance, if the $L_i$
  are the vanishing cycles of a K{\"a}hler degeneration and satisfy a
  relation as in (\ref{dickisgay}) the result reduces to
  an old theorem of Werner \cite{W}.  \medskip

\noindent Here is an illustrative (and suggestive) example.  
Cohomology classes of K\"ahler forms always have special properties
  \cite{GH}.  For instance, if $\omega$ is a K\"ahler form on a
  six-manifold then the 
pairing $H_4{\times}H_4\rightarrow\R$ given by $(A,B) \mapsto
A \cap B \cap \mathrm{PD}[\omega]$ is non-degenerate (Hard Lefschetz
theorem) and of signature $(1+2h^{2,0}, h^{1,1}-1)$ (Hodge-Riemann
bilinear relations).  If $b_2(X)=3$ the latter condition implies that
the matrix of the bilinear form $\cap [\omega]$ has positive
determinant, having an odd number of positive eigenvalues.
Fix a split K\"ahler
  structure $\b\omega_{\P^2} \oplus \a\omega_{\P^1}$ on
  $\P^2\times\P^1$, with $\a,\b >0$.  The manifold
  contains a Lagrangian sphere $L$, the graph of the composition of the Hopf
  map and the antipodal map on $\P^1$, precisely when $\b>\a$.  (To see
  this, note that the volume $\b^2\pi^2/2$ of $\P^2$ is exactly filled
  by a symplectic ball of radius $\sqrt{\b}$ in $\C^2$, whilst the
  boundary of the ball of radius $\sqrt{\a}$ 
  induces a form on $\P^1$ of volume $\pi\a$.)  Both conifold
  transitions along $L$ admit symplectic structures; if we flop the
  resolving sphere then change $D \mapsto -D$, in the notation of the
  proof of Theorem (\ref{dickisgay}).  The symplectic forms guaranteed
  by the theorem
  have the form $\Omega_{\varepsilon} = \pi^*\omega \pm \varepsilon
  \mathrm{PD}[D]$, for some small $\varepsilon$.  By computing Chern
  numbers, one can see that the conifold transitions are not given by
  (say) blowing up a rational curve inside $\P^1\times\P^2$, so they seem to
  have no obvious K\"ahler construction.  
 
\begin{prop}\label{salvage}
The determinant of the matrix given by cap product with
$\Omega_{\varepsilon}$ is positive for small
$\varepsilon$ if and only if $\b>\a$.
\end{prop}

\begin{proof}
We will work in the standard complex orientations on both
factors.  However, to simplify some formulae and remove the antipodal
maps, we will take the flipped symplectic structure $\a\omega_{\P^1}
\oplus -\b\omega_{\P^2}$.
To describe the intersection form of the conifold
transitions, we will first need to describe the cycles on these
spaces.  According to \cite{At}, the twistor
space of $\R^4$ is isomorphic to the total space of $\mathcal{O}(1)
\oplus \mathcal{O}(1) \rightarrow \P^1$, so taking duals we obtain an
isomorphism between $\mathcal{O}(-1)^{\oplus 2}$ and
$\SS^2_J\times(\R^4)^*$.  Here the first factor parametrises the
complex structures on $(\R^4, \sum dx_i \otimes dx_i)$ and at a point
$j\in\SS^2_J$ the complex structure on $\{j\} \times (\R^4)^*$ is that
induced from the natural isomorphism (over $\R$) between the complex dual of
$(\R^4,j)$ and the real dual space $(\R^4)^*$.  Combining with
(\ref{cotangentisnode}) and bearing in mind the discussion of
(\cite{STY}, Appendix) we obtain a map
\begin{equation} \label{project}
\SS^2_J \times (\R^4\backslash\{0\}) \ \rightarrow \ T^*\SS^3;
\qquad (j,v) \mapsto (v/|v|,jv).
\end{equation}
This is a diffeomorphism off zero-sections, and we use this map to
exhibit a small resolution; in other words, we replace $L$ by
$\SS^2_J$ via this map. \medskip

\noindent Consider the 3-sphere $U = \{([w],w) \in \P^1 \times \SS^3 \}$ where
$\SS^3 = \partial B\subset
\P^2$ is the boundary of a standard embedded Euclidean ball $B^4$ of
appropriate radius.  Fix a vector $v \in \C^2$ of length 1.  $U$
intersects $[v]\times\P^2$ transversely along a Hopf circle $\SS^1 = \{([v],
e^{i\theta v})\}$, and locally about this $\SS^1$ the copy of $\P^2$
can be smoothly identified with 
$$T^*\SS^3|_{\SS^1} \, \cong \, \big\{ (e^{i\theta v}, \lambda j
  e^{i\theta v})\, | \, \lambda \in 
  (0,\infty), j \in \SS^2_J \big\} \subset T^*\SS^3.$$
Via the map (\ref{project}) this corresponds to a (non-complex) real
  rank two bundle 
\begin{equation} \label{localnearsphere}
\SS^2_J \times \R\langle v,iv\rangle \subset
  \SS^2_J \times \R^4.
\end{equation} 
Let $|\P^2|$ denote the proper transform of the chain $[v]\times\P^2$ in
the small resolution, given by taking the closure of its image under
(\ref{project}) when we include back the zero-section on the LHS.
Since $[v]\times\P^2$ and $[w]\times\P^2$ are disjoint for
$[v]\neq[w]$ the proper transforms meet only along the exceptional
curve $\SS^2_J$, and hence the triple intersection $(|\P^2|)^3$ is
given by the Euler class of the normal bundle of $\SS^2_J$
inside $|\P^2|$.  However, we have argued that locally near the
resolving sphere $|\P^2|$ looks like (\ref{localnearsphere}), and
  hence this Euler class is trivial: $|\P^2|^3=0$. \medskip

\noindent Now let $R$ denote the chain $\{([w],w) \in \P^1\times\C^2 \, | \,
|w|\leq 1\}$ with boundary $U$.  Abstractly, this is isomorphic to the
disc bundle of the $\mathcal{O}(-1)$ line bundle over $\P^1$, or the
blow-up of the disc $D^4$ at the origin.  Write $|R|$ for its proper
transform in the small resolution.
We claim that $|R|^2 \cdot |\P^2| = -1$.  First, we
try to move $R$ off itself.  Let $\Phi_t$ be the flow $w \mapsto
w\cos(t)+jw\sin(t)$ giving 
$$R' = \big\{ ([w],\Phi_{1-|w|}(w)) \in \P^1 \times B^4 \, | \,
|w|\leq 1 \big\}.$$
This flow fixes the $\P^1$ at the origin, the $U=\SS^3$ at the
boundary, and nothing else.  Hence $|R'|\cap|R|$ comprises the $\P^1$
at the origin, together with a contribution from $\SS^2_J$;
moreover, locally about the $\P^1$ both $|R|$ and $|R'|$ look like
$\mathcal{O}(-1)$.  We have already observed that $|\P^2|\cdot
\SS^2_J = 0$ and hence $|R|^2\cdot |\P^2|$ can be computed inside
$\P^2 \times \P^1$ as $([v]\times\P^2) \cdot [-\P^1] = -1$.  To
justify the sign, argue as follows.  If $[w] \in \P^1$ is fixed, the
cycle $R$ is given by a fibre $\langle w,iw \rangle$ of
$\mathcal{O}(-1)$, which is a complex line, whilst $\Phi$ is a
perturbation in the direction $\langle jw,jiw=-kw \rangle$, where
these span an anticomplex line.  It follows that $R$ and $R'$ meet
transversely but negatively along the zero-section. \medskip

\noindent From these two intersection calculations, we can deduce the
result.  Take as ordered basis for $H_4$ of the conifold transition the cycles
$|\P^2|, |\P^1\times\P^1|$ and $|R|$.  Thinking of the
$\P^1\times\P^1 \subset \P^1\times\P^2$ as having a line at
infinity in the second factor, it is immediate that
$|R| \cap |\P^1\times\P^1| = \emptyset$, from which we deduce:
$$|\P^1\times\P^1|\cdot|\P^2|^2 = 0; \
|\P^1\times\P^1|^2\cdot|\P^2|=1; \
|\P^1\times\P^1|\cdot|R|^2=|\P^1\times\P^1|^2\cdot|R|=0.$$
Denote the remaining triple intersection numbers
$|\P^2|^2\cdot|R|=A$ and $|R|^3=B$ for integers $A,B$.  From (\cite{STY}, Proof
of Theorem 2.9) it follows that $A=\pm 1$, since near $U$ our four-chain $R$ is
exactly a collar neighbourhood as described in (op.cit).  (Locally,
then, $R$ defines a complex surface in the small resolution which
either meets the exceptional $\P^1$ transversely in a point or
contains it with normal bundle $\mathcal{O}(-1)$.)  The
symplectic forms on the two conifold transitions have Poincar\'e dual
homology classes 
$$\mathrm{PD}[\Omega_{\varepsilon}] \ = \
\a[|\P^2|]-\b[|\P^1\times\P^1|] \pm \varepsilon[|R|]$$
with $\a,\b$ both strictly positive.  (Note that here the second
term is $-\b$ because with respect to $-\omega_{\P^2}$ we reverse
the sign of the generator of $H_2(\P^2)$.  However, the only
non-trivial triple product involving the four-cycle $\P^1\times\P^1$
is $|\P^1\times\P^1|^2 \cdot |\P^2|$, where it occurs as a square
and the sign is irrelevant.)  In the given ordered basis, the
matrix which encodes the bilinear form $\cap[\Omega_{\varepsilon}]$ on
$H_4$ is the following:
$$
\left(
\begin{array}{ccc}
\pm \varepsilon A & -\b & \a A \mp
  \varepsilon \\
-\b & \a & 0 \\
\a A \mp \varepsilon & 0 & -\a \pm B\varepsilon 
\end{array}
\right) \qquad = \qquad 
\left(
\begin{array}{ccc}
0 & -\b & \a A \\
-\b & \a & 0 \\
\a A & 0 & -\a
\end{array} \right) \ + O(\varepsilon). 
$$
The determinant of this matrix is $\b^2\a - \a^3A^2 +
O(\varepsilon)=\a(\b^2-\a^2)+O(\varepsilon)$, which is positive if $\b > \a$, and this
completes the proof. 
\end{proof}

\noindent The possibility of finding non-K\"ahler examples of the
surgery is made more interesting by 
a question of Donaldson \cite{Do}, who asked if every
Lagrangian sphere in a complex algebraic variety arises as the
vanishing cycle for some K\"ahler degeneration.  Given a symplectic  
six-manifold $X$ containing a Lagrangian
sphere which is trivial in homology, we can form the conifold 
transition and the blow-up of the nodal space, and by Theorem
\ref{dickisgay}, resp. 
Section 2, each has a distinguished deformation class of 
symplectic forms (at least once we fix the four-chain bounding the
Lagrangian).  

\begin{prop}
In the situation above, suppose that $X$ is K\"ahler with
$h^{2,0}=0$, and that the Lagrangian can be collapsed in a K\"ahler degeneration.
Then the above symplectic forms on the conifold transition are K\"ahler. 
\end{prop}

\begin{proof}
$h^{2,0}$ is constant in the K\"ahler degeneration (by the constancy of
$h^{2,0}+h^{1,1}+h^{0,2}$ and the upper semicontinuity of each term), and by
standard technology (e.g. mixed Hodge structures)
$h^{2,0}$ remains zero on the big resolution. Also using $h^{2,0}=0$ we see that
$X$ and the big resolution are in fact \emph{projective},
with hyperplane class $[H]$ arbitrarily close to a high multiple of any given
symplectic form in the distinguished family. So we may as well choose the
original symplectic form on $X$ to be projective.

Since $[L]=0 \in H_3(X)$ the
proper transform of a chosen bounding four-chain gives a cycle on the
big resolution which hits one $\P^1$ in the exceptional divisor and
not the other (\cite{STY} Theorem 2.0); $h^{2,0}=0$ implies that this cycle can be
taken to be a divisor $D$. 

Since $D$ intersects one ruling of the exceptional divisor but on the other,
we can find a linear combination of $D$ and $H$ which is ample except that it
evaluates to zero on precisely one of the two rulings. Taking the map to $\P^N$
associated to the linear system of a high multiple of such a class
will contract (only) this ruling and yield the small resolution. By varying
the linear combination we get all (high multiples of) the symplectic
forms $\omega_X + \varepsilon \sigma$ (using the notation of the proof of
Theorem \ref{dickisgay}) with $\varepsilon$ rational. The others follow
by continuity.
\end{proof}

\noindent Hence, if the (distinguished deformation class of forms on
the) conifold transition of a
projective variety with $h^{2,0}=0$ (e.g. $\P^1\times\P^2$) in a
single Lagrangian $L$ is non-K{\"a}hler, this
Lagrangian is not a vanishing cycle for any K\"ahler degeneration.
(More general conifold transitions presumably give rise to configurations of
Lagrangians that are not simultaneously realised as vanishing cycles.)
This opens another approach to proving certain Lagrangian spheres are
not vanishing cycles of K\"ahler degenerations:  employ Lemma
\ref{perhapsfano} to see that the full blow-up is Fano, and then
invoke the Mori-Mukai classification of Fano 3-folds. (We would like to
use Fano-ness since there is no general classification of
3-folds and, as above, it seems hard to violate other topological
constraints such as the Hard Lefschetz theorem
via conifold transitions.)
  
The following is also relevant to Donaldson's
question, tackling it from a Calabi-Yau, rather than Fano, perspective. 
Recall that a Calabi-Yau 3-fold is \emph{rigid} if $b_3=2$;
equivalently $h^{2,1}=0=h^1(TX)$, i.e. it has no complex deformations.

\begin{lem} \label{wontvanish}
An essential Lagrangian sphere in a rigid Calabi-Yau 3-fold is never a
vanishing cycle.
\end{lem}

\begin{proof}
If we have a 
flat family of varieties over the disc, with one fibre $Z$ and central
fibre nodal, then we may pick a fibrewise holomorphic
3-form varying holomorphically and not tending to zero at the central
fibre (the pushdown of the relative canonical sheaf is torsion free
rank 1, and so a line bundle, so we may pick a nowhere zero section).
Thus the period $\int_L\Omega$ does not tend to zero and $L$ is not
collapsed. 
\end{proof}

\noindent If we could find sufficiently ``degenerate''
Lagrangian spheres, we would certainly have non-algebraic surgeries.

\begin{lem}
If $L \subset X$ is a Lagrangian three-sphere which bounds a smoothly
embedded four-ball in $X$, the conifold transition of $X$ in $L$ is
not homotopy K{\"a}hler.
\end{lem}

\begin{proof}
On collapsing the $\SS^3$, the $D^4$ becomes an $\SS^4$ that lifts to one of
the small resolutions \cite{STY} (and $\overline{\P^2}$ in the other).
Its $H_4$ class is nonzero (it intersects the exceptional $\P^1$ in $+1$) but
has intersection zero with all other $H_4$ classes (as the intersection factors
through $H_2(\SS^4)=0$).  This would contradict the Hard Lefschetz theorem,
so the manifold is not K{\"a}hler, though it is symplectic by
(\ref{dickisgay}).
\end{proof}

\noindent It is not clear if such ``contractible'' Lagrangian
3-spheres can ever exist in closed symplectic six-manifolds, although
the existence of such spheres in fake symplectic $\R^6$s \cite{Mu}
indicates that any obstruction would be global.
In the Calabi-Yau context, there are no known
examples even of homologically trivial Lagrangian 3-spheres. 
Nonetheless, the existence of the symplectic conifold
transition makes it very plausible that there are families of
symplectic non-K{\"a}hler Calabi-Yaus.  For instance, the quintic
hypersurface $Q \subset \P^4$ is a Calabi-Yau with $b_3(Q)=204$.
There are well-known examples of nodal quintics containing as many as
130 nodes \cite{vS}, whose projective small resolutions are
rigid.  

\begin{ques}
Is there a rigid Calabi-Yau 3-fold $Q$ containing a Lagrangian
3-sphere ?
\end{ques}

\noindent As in the previous section, all answers to the question are
interesting.  Suppose first the answer is positive, and that there are
two disjoint essential Lagrangian
spheres which (necessarily) satisfy a relation as in
(\ref{dickisgay}).  Then the small resolution 
given by the Theorem above would be a symplectic
Calabi-Yau 3-fold with \emph{trivial} third Betti number.  This could not be
K{\"a}hler, since the cohomology class of a holomorphic volume form
on a Calabi-Yau is necessarily non-zero.  
If an essential sphere exists but
satisfies no good homology relation, we obtain examples of 
Lagrangian spheres in algebraic varieties that can never arise as 
vanishing cycles for complex degenerations, via (\ref{wontvanish}).
Finally, as remarked above, there are no known examples of
homologically trivial 3-spheres in Calabi-Yaus, so even the existence
of an inessential sphere in $Q$ would be novel.  

In a slightly
different direction, one can show that the existence of the surgery
(\ref{dickisgay}) implies that either there are quintic 3-folds
with $j$ nodes for every $j \leq 130$ -- itself rather surprising --
or there are non-K{\"a}hler symplectic Calabi-Yaus, given by
conifold transitions on the examples from \cite{vS}. \bigskip

\section{Fibre products and triple points}

Let us consider one more complicated instance of passing
from a configuration of Lagrangian vanishing cycles to a new
symplectic manifold containing a symplectic exceptional divisor.  This
will show that these more complicated surgeries can indeed be amenable
to explicit 
computation and construction.  
Again motivated by a desire to find new symplectic manifolds with
$c_1=0$, we pass from double points to 3-fold \emph{triple
  points}, that is isolated singularities of the form 
$$R = \{ x^3 + y^3 + z^3 + w^3 = 0 \}.$$
The singularity at the origin can again be resolved by a single
blow-up, and now the exceptional divisor $E$ is a cubic surface in $\P^3$
(abstractly diffeomorphic to the six-fold blow up of $\P^2$, so
$b_2(E)=7$). Its normal bundle in the resolution is $\mathcal{O}_{\P^3}(-1)|_E=
K_E$, so by the adjunction formula the blow-up has trivial canonical
bundle over $E$;
in particular the transition preserves the Calabi-Yau condition.

Before describing the smoothing of the singularity, it will be helpful
  to introduce a pretty construction of
  Lagrangian spheres (hence of degenerations, or Lagrangian
  blow-downs), often appropriate to the Calabi-Yau setting. 
Given a pair of smooth surfaces $S_i$ fibred over a curve $C$, we can form
their fibre product $S_1\times^{\ }_C S_2$, as used so effectively by
  \cite{Sch2}, 
for example. To analyse its singularities, we look at the local model
$$
S_i=\{f_i(x_i,y_i)=t)\}\subset \C^2\times \C,
$$
fibred over $\C$ by the $t$ variable. The fibre product, then, is locally
the threefold
$$
S_1\times^{\ }_C S_2=\{f_1(x_1,y_1)=f_2(x_2,y_2)\}\subset\C^2\times\C^2.
$$
So we see that there are only singularities $(x_1,y_1,x_2,y_2)$ if \emph{both}
points $(x_i,y_i)$ lie on singular fibres above the \emph{same} point $t$
of $C$, at the singular points. So for instance if both $f_i$ define double
point singularities in the curve fibres over $t=0$, $f_i(x_i,y_i)=x_iy_i$, say,
then the threefold also has a double point $x_1y_1=x_2y_2$. More generally
degree $n$ curve singularities $f_i=x_i^n+y_i^n$ give degree $n$ threefold
singularities $x_1^n+y_1^n=x_2^n+y_2^n$. We discuss the $n=3$ triple point
case presently.

\begin{figure}[h]
\center{
\input{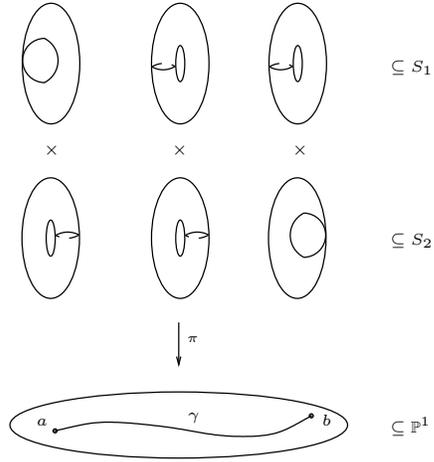}
\caption{Fibred Lagrangian three-spheres. \label{S3}}}
\end{figure}

To smooth such singularities we can first Morsify the $f_i$ to reduce to
double points, then move one fibration by an automorphism of the base to
move its singular fibres away from those of the other surface. The Lagrangian
$\SS^3$ vanishing cycle of this latter smoothing is easily
described. Take a path 
$\gamma$ in the base from the image of the singular fibre of $S_1$ to that
of $S_2$. Over this, by symplectic parallel transport, lies a fibration
by $\SS^1$ vanishing cycles for the curve double point in the curve fibres
of $S_i$. Taking the fibre product of these two fibrations over $\gamma$
gives a $\rot^2=\SS^1\times\SS^1$-fibration with the property that one
$\SS^1$ factor collapses at one end of $\gamma$ and the other at the other
end. See Figure \ref{S3}. Thus over each half of $\gamma$ we get a
solid torus, glued together
to form an $\SS^3$ handlebody. This is easily seen to be Lagrangian by
putting the fibre product into the full product, equipped with the
product symplectic form.

\begin{figure}[h]
\center{
\input{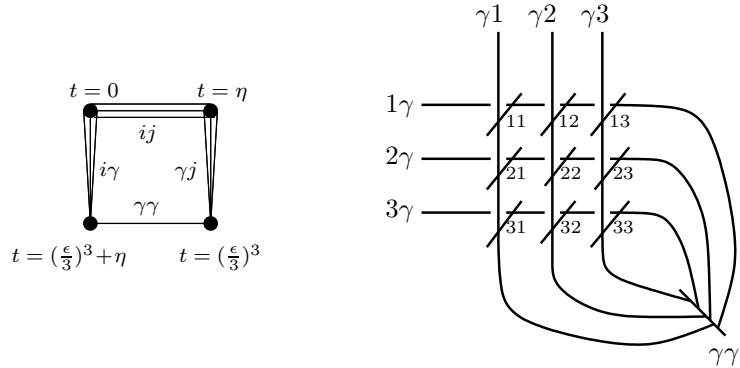}
\caption{The vanishing cycles of the 3-fold triple point, $\rot^2$-fibred over the
paths in $\P^1$ on the left \label{triple}}}
\end{figure}

With this background, we can study the smoothing of the 3-fold triple point
singularity.  A helpful picture of the configuration of vanishing
cycles, similar to a description given already by Ebeling \cite{Eb}, is given
in the following (note that the singularity has Milnor number 16):

\begin{lem} \label{triplepoints}
Let $R'$ be the smoothing of a triple point singularity. Then $R'$
contains a configuration of 16 Lagrangian spheres with intersections
as indicated in Figure \ref{triple}.
\end{lem}

\begin{proof}
We take as model for the triple point a fibre product of two elliptic
fibrations $E_i \rightarrow \P^1$ each of which contain a
two-dimensional triple point $\{xy(x+y)=0\}$ inside a torus fibre.
So locally the surfaces look like $xy(x+y)=t$ and $uv(u+v)=t$; the threefold
fibre product then is locally $xy(x+y)=uv(u+v)$ -- a triple point.

To smooth this, we move the fibre products apart by changing $t$ to $t-\eta$
in the second fibration, and then Morsify them both by replacing $xy(x+y)=t$
by $xy(x+y+\varepsilon)=t$. (So the threefold is now locally $xy(x+y+\varepsilon)=
uv(u+v+\varepsilon)+\eta$.) This new surface has only nodes in its
fibres, at
$$
(x,y,t)\in\left\{(0,0,0),\,(0,-\varepsilon,0),\,(-\varepsilon,0,0),\,\left(
-{\varepsilon\over3},-{\varepsilon\over3},\big({\varepsilon\over3}\big)^3\right)\right\}.  
$$
That is, the fibre
over $t=0$ is 3 $\P^1$s, with the vanishing cycle of the node over
$t=\big({\varepsilon\over3}\big)^3$ being the essential loop (``triangle")
in this triangle of $\P^1$s; this is contracted to the original curve
triple point (3 coincident lines) as $\varepsilon\to0$.  (In Kodaira's
notation \cite{BPV}, these triangles of rational curves are called
$I_3$ singular fibres.)

Labelling the 3 nodes over $t=0$ (and their $\SS^1$ vanishing cycles
in nearby fibres) 
by $1,2,3$, and the other vanishing cycle by $\gamma$, we get 16 Lagrangian
$\SS^3$ vanishing cycles by taking any one of these four in the first
surface and 
any one of them in the second, and taking their fibre product to give a
$\rot^2$-fibration over the path between the points in $\P^1$ at which they
collapse, as above (Figure \ref{S3}).
Since $\gamma$ intersects each of $1,2,3$ in one point and there
are no other intersections in the fibres, the intersection pattern of the
resulting $\SS^3$ is easily determined to be as in Figure \ref{triple}.
\end{proof}

\noindent Now given two elliptic fibrations, we can often detect these
triple point degenerations by hand, and hence give examples of the
triple point transition.  One example is given by studying one of
Schoen's rigid Calabi-Yau fibre products from \cite{Sch2}.  There is a
rational elliptic surface $\pi: E(1) \rightarrow \P^1$ with four
$I_3$ singular fibres.  This can be obtained from the pencil of cubic
curves $\{ x^3 + y^3 + z^3 + 3txyz = 0 \}_{t \in \P^1}$, which has an
$I_3$ fibre when $t\in\{e^{i\pi/3},e^{5i\pi/3},-1,\infty\}$.
Let $Z$ denote the fibre product
$E(1) \times_{\pi} E(1)$, which has 3-fold nodes at the 36 points
$(N,N')$ with $N, N'$ nodes in some fibre $\pi^{-1}(t)$ of $\pi$.
There is a small resolution $\tilde{Z}$ of $Z$ at these 36 points which is a 
projective rigid Calabi-Yau 3-fold. 

\begin{lem}
In a suitable basis, the monodromy of the elliptic surface $E(1)$ is given by 
$$t(a)^3 \cdot t(b)^3 \cdot t(a+b)^3 \cdot t(2a+b)^3 = 1$$
where $a,b$ are the standard meridian / longitude curves on $\rot^2$ and
$t(C)$ denotes the Dehn twist in a simple closed curve in the homology
class $C$.  (Here the generating loops are ordered clockwise around
the base-point.)
\end{lem}

\begin{proof} 
It is well-known that one
can construct a holomorphic elliptic surface from any appropriate
monodromy representation, and moreover that there is a \emph{unique}
elliptic surface with exactly four $I_3$ singular fibres (see
\cite{Mi} and \cite{Pe}). Hence it is enough
to see that the product of the four cubes of Dehn twists above is
indeed the identity.  In standard conventions one has
$$
t(a) = \left( \begin{array}{cc} 1&1 \\ 0&1 \end{array} \right); \ 
t(b) = \left( \begin{array}{cc} 1&0 \\ -1&1 \end{array} \right); 
$$
$$
t(a+b) = \left( \begin{array}{cc} 0&1 \\ -1&2 \end{array} \right); \
t(2a+b) = \left( \begin{array}{cc} -1&4 \\ -1&3 \end{array} \right)
$$
from which the result follows by a direct computation.
\end{proof}

\noindent It is important to notice that the $a$ and $b$ vanishing
cycles are \emph{adjacent} in the monodromy representation; if we perturb the
$I_3(b)$ fibre into an $I_1$ and an $I_2$ fibre, the $I_3(a)$ and the
$I_1(b)$ give the right configuration of cycles to degenerate to an
elliptic surface with a fibre containing a triple point, as in the
proof of Lemma (\ref{triplepoints}).  \medskip

\noindent This perturbation is certainly possible
symplectically, so separate the $I_3(b)$ singular fibre into an $I_1$
and $I_2$.  We can re-degenerate $E(1)$ to an elliptic surface with a
triple point singularity.  The new fibre product has a 3-fold triple
point, which can be blown up preserving $c_1=0$.  In this
case, the resulting manifold is apparently K{\"a}hler; the degenerations
can all be realised holomorphically, since monodromy data
is enough to 
guarantee the existence of rational elliptic surfaces, according to
the theory of Miranda and Persson \cite{Mi,Pe}.  However, this is not to
say that the surgery is uninteresting.  The manifold we obtain is a
new rigid Calabi-Yau which is not itself the small resolution of a
fibre product 
of elliptic fibrations, and hence not on Schoen's list.  (Blowing up a
triple point singularity increases the Euler characteristic by 24,
from which one can compute its Betti numbers.)  In general
constructions of rigid K{\"a}hler Calabi-Yaus are not plentiful, and it
is likely that other examples can be found by arguments such
as above.  In any case, the fibred nature of these three-folds makes
it very easy to find large collections of Lagrangian spheres (and
compute their intersections from intersections of plane curves),
something that is of considerable importance in the general programme
of mirror symmetry. \bigskip

\section{A four-dimensional interlude}

So far we have always passed from a smoothing to a resolution.  If we
go in the other direction, we can try (as in the minimal model programme) 
to eliminate curves on which the canonical class is negative, and 
find some consequences for symplectic manifolds
which are neither Fano nor Calabi-Yau but of ``general type''; that
is, for which the canonical class of a given almost complex structure
itself contains symplectic forms.  In general, the difficulty here is
that one needs a family of symplectic structures which degenerate (as
forms, and not just in volume) along the locus which one wants to
contract; as yet, there is no good theory of ``symplectic extremal
rays'' in this general sense.  However, this situation is often
provided by algebro-geometric arguments.  The canonical class of a
complex surface cannot contain a K{\"a}hler form if the surface contains
a holomorphic $-2$-curve.  By contrast:

\begin{prop}
If $X$ is a complex surface of general type, then $X$ has symplectic
general type.
\end{prop}

\begin{proof}[Sketch]
The multicanonical linear system is an embedding away from ADE trees
of rational curves which are contracted to isolated singularities.
These have local smoothings in which the complex curves become
Lagrangian two-spheres.  By replacing
one model with the other, in the vein of the first section, and
noticing that the symplectic form is cohomologically unchanged since
exact in a neighbourhood of the vanishing cycles, one
quickly arrives at the result above.  For details see 
\cite{STY} (the result was independently proven by Catanese
in \cite{Ca}, although the former proof is more relevant for the
discussion here).  
\end{proof}

\noindent The small resolution of 3-fold double points goes over to the
\emph{simultaneous resolution} of surface double points, due to
Brieskorn.  This has been
notably exploited by Seidel \cite{Se1} and Kronheimer \cite{Kr}.  There is
a family of symplectic manifolds over the disc $D$ which for every
$t\neq0$ contains a Lagrangian two-sphere, and which at $t=0$ contains
a rational $-2$-curve (on which the symplectic structure is completely
degenerate) -- the resolution of the singular point fits
into a family of smoothings (after passing to a double cover of the
base). The total space of this family, containing the rational curve
in the central fibre, is exactly the small resolution of a 3-fold
ordinary double point which arises from the base change
$\{x^2+y^2+z^2=t^2\} \mapsto \{x^2+y^2+z^2=t\}$.  The smooth monodromy of
the family $\{ \sum_{j=1}^3 x_j^2 = \varepsilon \}$ over the disc 
has order two, but the symplectic monodromy -- given by a generalised
Dehn twist in the Lagrangian vanishing cycle -- has infinite order
\cite{Se1}.  \medskip

\noindent The resolution and smoothing are diffeomorphic precisely for
simple singularities, so if we move beyond these then we obtain
surgeries which have a non-trivial topological effect.  For instance, 
there are examples analogous to that above in which the
isolated complex curve has higher genus:  a degree $d$
surface singularity has a resolution with a degree $d$ curve in
$\P^2$ as exceptional set.  This leads to a surgery in which
configurations of Lagrangian spheres in symplectic four-manifolds can
be blown down and replaced by symplectic surfaces of high genus and
negative square.  In particular cases, these surgeries (or more
properly their inverses) can be related to familiar 
operations on symplectic four-manifolds.  We give the following
proposition in the case $d=3$, parallel to the discussion of triple
points in 3-folds in the previous section, but it is not hard to
generalise to arbitrary degree.  \medskip

\noindent The
resolution of the surface triple point $\{x_1^3+x_2^3+x_3^3=0\}$ is a
genus one curve with 
normal bundle having first Chern class $-3$.  Introduce the notation
$\mathcal{C}_3$ for an open neighbourhood of the configuration of vanishing
cycles of the smoothing, which we can assume is diffeomorphic to an affine
cubic.  The singularity has Milnor number 8, but --
as with all non-simple surface singularities -- the
intersection matrix of the vanishing cycles is not negative definite.
Recall also that we can define the \emph{proper transform} of a
symplectic surface $C^2 \subset X^4$ inside the blow-up of $X$ along
$C$ by making $C$ 
$J$-holomorphic for some compatible almost complex structure $J$ on
$X$ integrable near $C$, and then taking the usual K\"ahler blow-up
and the holomorphic proper transform.

\begin{prop}
Let $X$ be a symplectic four-manifold which contains a symplectic
torus $\rot^2$ of square 
zero.  The following two operations are smoothly equivalent.  (i) Blow
up $\rot^2$ three times, and then replace a tubular neighbourhood of
its proper transform (a square $-3$ torus) by
$\mathcal{C}_3$. (ii) Fibre sum $X$ with a rational elliptic
surface $E(1)$ along $\rot^2$ and a fibre respectively.
\end{prop}

\begin{proof}
We regard the total space of the smoothing $\{x_1^3+x_2^3+x_3^3=t\}$, for an
appropriate polynomial $t$, as the complement of the
hyperplane at infinity in a cubic surface $Z=\{\sum_{j=1}^3 x_j^3 =
tx_4^3\} \subset \P^3$.  Topologically $Z$ is the
six-fold blow-up of 
$\P^2$, and the hypersurface 
$\{x^3+y^3+z^3=0\}$ is a torus of square $3$.  The complement
$Z\backslash\rot^2$ of the torus -- which has $b_2=8$ -- 
can be identified with the complement of a fibre and three sections in
the rational elliptic surface $\P^2{\#}9\overline{\P}^2 \cong
E(1)\rightarrow \P^1$.  The result now follows from the following
equivalence:  fibre summing the square $3$ torus given by blowing up six
base-points in a pencil of cubics with a $-3$ torus given by blowing
up a square zero torus three times, is the same as fibre summing $E(1)$
along the original square zero torus.  This equality in turn follows
from considering Gompf's pairwise version of the symplectic sum, which
enables one to perform the first fibre sum so the three $-1$-curves transverse
to the $-3$ torus are glued onto the three base points of 
the linear system of square three tori.
\end{proof}

\noindent Thus, in the situation of the Proposition, the
``degenerate-resolve'' surgery is smoothly equivalent to
``de-fibre-summing'' with a 
copy of $E(1)$ and then blowing up three times.  
If $d>3$ there is a similar interpretation; the
``contract-deform'' direction of the surgery -- contracting a curve $C
\subset X$ of genus $g=(d-1)(d-2)/2$ and square $-d$ which was obtained
as a $d$-fold blow-up of a square zero curve -- is smoothly equivalent to
a fibre sum of $X$ along a particular
genus $g$ Lefschetz fibration.  The fibre-sum with $E(1)$ has been an
extremely useful surgery \cite{FS} and it is encouraging that it can
be recovered from this point of view.  \bigskip

\section{Tree-like configurations}

After isolated Lagrangians, the simplest situation to 
consider is that of contracting linear chains of spheres.  In four real
dimensions the $A_n$ singularities have diffeomorphic resolutions and
smoothings, but this is no longer true in six dimensions.      
These surgeries, however, still don't produce any more 
(diffeomorphism types of) symplectic manifolds than those one can obtain from
the ordinary double point surgery, at least when working with small resolutions.  One
can see this as follows.  

\begin{lem} \label{link}
The link of an $A_n$ chain of spheres is diffeomorphic to $\SS^5$ if
$n$ is even and $\SS^2 \times \SS^3$ if $n$ is odd.
\end{lem}

\begin{proof}
The link is a series of connect sums of $\SS^3\times \SS^2$s, across
$\SS^2\times\SS^2$s (given by removing a ball in $\SS^3$ to give
$\SS^2$ boundary, and timesing
everything by $\SS^2$). So inductively it is enough to prove that
$$
(\SS^3\times \SS^2)\#^{\ }_{\SS^2\times \SS^2}(\SS^2\times\SS^3)\,
\cong\, \SS^5.
$$ 
Here the ordering of the factors is meant to indicate that the first
factor in the gluing locus $\SS^2\times
\SS^2$ is the boundary of a ball in the first $\SS^3$ factor,
but is the full $\SS^2$ factor in the second $\SS^3\times \SS^2$,
and the opposite 
for the second $\SS^2$ factor in the gluing locus. In the same
notation, the connect sum is 
$$
(D^3\times \SS^2)\cup_{\SS^2\times\SS^2}^{\ }(\SS^2\times D^3),
$$
since $\SS^3$ minus a ball is $D^3$. But this is isomorphic to taking
$\SS^2\times\SS^2\subset\R^5\subset\SS^5$ and filling it inside $\R^5$
to give the first $D^3\times \SS^2$ in the union, and filling it `outside'
in the $\SS^5$ to give the other.  (This is analogous to realising
$\SS^3$ as a union of two genus one handlebodies.  The embedding
$\SS^2\times\SS^2 \subset \SS^5$ comes from $\{(z_1,z_2,z_3)\in\C^3 \,
| \, |z_1|^2+|z_2|^2 = \frac{1}{2} = |z_2|^2+|z_3|^2\}$.)
\end{proof}

\noindent We employ this as follows.  It is well-known that the
singularity given by contracting an $A_n$-chain of Lagrangian
3-spheres has a small resolution if and only if $n=2k+1$ is odd.
This resolution is given by a $k$-fold cover of the ordinary double
point resolution, lifting the obvious cover $\{ \sum x_i^2 = t^{2k} \}
\rightarrow \{ \sum x_i^2 = t^2 \}$ given by $t \mapsto t^k$.  The
exceptional locus is a $k$-times thickened $\P^1$.   In this 
case, we can talk about conifold transitions in the entire $A_n$-chain, but:

\begin{lem} \label{boundary}
The conifold transition of $X$ in the $A_n$-chain ($n=2k+1$ odd) is
diffeomorphic to a manifold obtained by performing conifold
transitions in isolated Lagrangian spheres.
\end{lem}

\begin{proof} 
Algebraically, we can collapse alternate Lagrangians in
the chain by using the partial smoothing $\sum_{i=1}^3x_i^2=p_k(t^2)$ of
the full degeneration $\sum x_i^2=t^{2k}$; here $p_k$ is a degree $k$
polynomial with simple zeros.  Taking conifold transitions
along alternate spheres in the
$A_n$-chain is equivalent to taking small resolutions of the double
points in this partial smoothing.  The interpolating 3-spheres acquire
boundary, lifting to 
$\SS^2{\times}[0,1]$ homotopies between the exceptional loci of the
small resolutions.  The small resolution of the $A_n$ singularity
described earlier
arises in the limit  $p_k(t^2)\to t^{2k}$, and the $k$-times
thickened $\P^1$ is then the limit
of bringing the $k$ $\P^1$s from the ODP conifold transitions
together.
 
Alternatively, the conifold transition in the $A_n$-chain is given by
dividing out a suitable circle
action (determined by $k$) on the $\SS^2\times\SS^3$ boundary of
(\ref{link}), and then 
collapsing one factor of the resulting $\P^1\times\P^1$.  
This clearly gives a filling of the link 
$\SS^2\times\SS^3$ with one-dimensional $H_2$.
\end{proof}

\noindent As Eliashberg pointed out in \cite{El}, these chains and trees
of Lagrangians are nonetheless not without interest.  The contact
structures induced 
on $\SS^5$ and $\SS^2\times\SS^3$ are, for instance, distinct as
$n$ varies.  By looking at boundaries of more complicated trees, one
can obtain distinct contact structures on 5-manifolds -- distinguished by the
topologies of their fillings (or by contact homology?) -- in
cases where the underlying classical invariants of the plane fields
are equal. \medskip 

\noindent With the triple point above, we dealt with configurations of Lagrangian
spheres that arose from some projective degeneration, and the question of
the modulus of the symplectic structure in a neighbourhood of the
vanishing cycles did not
really enter into the question of whether a surgery exists.  For the
ADE contractions on surfaces, the K{\"a}hler form can be assumed
standard in a neighbourhood of the contracted spheres (\cite{STY}
using a result from \cite{Se3}). More
generally, one ingredient in running a symplectic surgery via a
transition in some collection of Lagrangians is the following
(well-known folk-theorem). 

\begin{prop} \label{maslov}
Let $\{ L_1, \ldots, L_n\}$ be a collection of Lagrangian spheres in
$X$ for which all intersections are transverse.  Suppose the
associated intersection 
graph is a tree (has no loops).  Then the symplectic structure in a
neighbourhood of the Lagrangians is unique to symplectomorphism.
\end{prop}

\begin{proof}
The proof is a ``plumbed'' version of Weinstein's theorem \cite{We}.
His argument shows that for any Lagrangian immersion $\phi:L \rightarrow X$
there is an immersion $\Phi$ from a neighbourhood of the zero-section in
$T^*L$ to $X$ such that $\Phi^*\omega_X = dp \wedge dq$.  To deduce
the result, recall that symplectomorphisms act transitively on pairs
of transverse Lagrangian subspaces of a symplectic vector space, so we
can assume at each 
intersection point (say of $L_i$ and $L_j$) $\Phi$ pulls back $L_i$ to a fibre in
the cotangent bundle of $L_j$ and vice-versa.  Fix an intersection point
$P$, giving two preimages $\phi^{-1}(P)$ lying in $L_i$ and $L_j$.  We
can define two box
neighbourhoods $(D^n \times D^n)_{i,j}$ of these preimages in
$\amalg_{j=1}^n T^*L_j$, 
the total spaces of embedded symplectic balls lying in fixed Darboux
charts, and such that
the first factor in the $i$-box describes 
fibres of $T^*L_i$ and the second factor in the $j$-box describes
fibres of $T^*L_j$.  By construction, $\Phi$ is assumed to be the
immersion (quotient map) which identifies $(D^n \times \{t\})_i$
with $(\{t\}\times D^n)_j$ (and is an embedding away from the boxes).
Given two forms $\omega_X, \tilde{\omega}_X$ both making the same $L_j$
Lagrangian, we obtain two such
models comprising a collection of copies of $T^*\SS^n$ with the standard
symplectic structure and a
quotienting relation describing the immersion on the boxes around
intersection points.  These models can be identified symplectically by lifting a
diffeomorphism $\amalg L_i \rightarrow \amalg L_i$ whose differential
matches the appropriate box neighbourhoods.  This descends to give a
diffeomorphism of the open neighbourhoods of $\phi(\amalg L_i)$ inside
$X$ which intertwines the two symplectic structures.
\end{proof}

\noindent At least for tree-like configurations,
this says that if a symplectic
manifold contains a graph of Lagrangian spheres that is
combinatorially the same as the graph obtained by smoothing
some complex singularity, then there is an associated symplectic surgery which
collapses the spheres and blows up the (locally analytic) singular
point that results.  There are well known lists and techniques for
identifying these configurations coming from singularity theory.  

In fact, it seems that a generalisation of the above Proposition is
not much harder.  For any configuration of Lagrangian spheres,
there is a homological invariant -- the relative intersection numbers
for a choice of orientations on all the Lagrangians -- but the argument
above suggests this is the unique invariant.  (Maslov classes play
no role, indeed are not defined, since locally the loops in a
configuration of Lagrangians are not locally spanned by discs.)  
Such more general configurations of Lagrangians can be produced from
simple configurations by Dehn twist 
automorphisms \emph{\`a la} Seidel \cite{Se1}. \\

\noindent {\bf Acknowledgements:} The authors are very grateful to Denis Auroux
for many helpful suggestions.


\end{document}